# ON THE TAYLOR COEFFICIENTS
# OF THE HURWITZ ZETA FUNCTION

Khristo N. Boyadzhiev

*Department of Mathematics, Ohio Northern University,*
*Ada, Ohio, 45810*
*k-boyadzhiev@onu.edu*

**Abstract**. We find a representation for the Maclaurin coefficients $\zeta_n(a)$ of the Hurwitz zeta-function $\zeta(s,a) = \sum_{n=0}^{\infty} \zeta_n(a) s^n$, $|s| < 1$, in terms of semi-convergent series

$$\zeta_n(a) = -1 + \sum_{k=n}^{\infty} (-1)^{k+1} \begin{bmatrix} k \\ n \end{bmatrix} \frac{B_{k+1}(a-1)}{(k+1)!},$$

where $B_n(x)$ are the Bernoulli polynomials and $\begin{bmatrix} k \\ n \end{bmatrix}$ are the (absolute) Stirling numbers of the first kind. When $a = 1$ this gives a representation for the coefficients of the Riemann zeta function. Our main instrument is a certain series transformation formula.

A similar result is proved also for the Maclaurin coefficients of the Lerch zeta function.

**1. Introduction. Exponential polynomials and the Exponential Transformation Formula for series.**

The exponential polynomials (or single variable Bell polynomials) $\phi_n$ can be defined by

$$\phi_n(x) = e^{-x}(xD)^n e^x, \quad n = 0, 1, \ldots \tag{1.1}$$

(where $(xD)f(x) = xf'(x)$). Equivalently

$$\phi_n(x) e^x = (xD)^n e^x = \sum_{k=0}^{\infty} \frac{k^n}{k!} x^k \tag{1.2}$$

One has

$$\phi_0(x) = 1, \ \phi_1(x) = x, \ \phi_2(x) = x^2 + x, \ \phi_3(x) = x^3 + 3x^2 + x, \text{ etc.} \tag{1.3}$$

These polynomials were first studied by S. Ramanujan (see [3, Chapter 3] and [4] for further details). All polynomials $\phi_n$ have positive integer coefficients, which are the Stirling numbers of the second kind $\{{n \atop k}\}$ (or $S(n,k)$), $0 \le k \le n$. Thus

$$\phi_n(x) = \sum_{k=0}^{n} \{{n \atop k}\} x^k \tag{1.4}$$

The polynomials $\phi_n$ form a basis in the linear space of all polynomials. One can solve for $x^k$ in (1.2) and write the standard basis in terms of the exponential polynomials:

$$1 = \phi_0, \ x = \phi_1, \ x^2 = -\phi_1 + \phi_2, \ x^3 = 2\phi_1 - 3\phi_2 + \phi_3, \text{ etc.} \tag{1.5}$$

If we set
$$x^n = \sum_{k=0}^{n} (-1)^{n-k} [{n \atop k}] \phi_k \tag{1.6}$$

then $[{n \atop k}] \ge 0$ are the absolute Stirling numbers of first kind. In particular,

$$[{k \atop 0}] = 0 \ (k > 0), \ [{k \atop 1}] = (k-1)!, \ [{k \atop k}] = 1. \tag{1.7}$$

More information on the Stirling numbers can be found in [7].

Suppose now that $f(z) = \sum_{n=0}^{\infty} a_n z^n$ is an entire function. Multiplying (1.2) by $a_n$ and summing for $n = 0, 1, \ldots$, one obtains the exponential transformation formula (ETF)

$$\sum_{k=0}^{\infty} \frac{f(k)}{k!} x^k = e^x \sum_{n=0}^{\infty} a_n \phi_n(x), \qquad (1.8)$$

(for details see [4]).

## 2. The Hurwitz zeta function.

The Hurvitz zeta function is defined for $\operatorname{Re} s > 1, a > 0$ by

$$\zeta(s, a) = \sum_{n=0}^{\infty} \frac{1}{(n+a)^s}. \qquad (2.1)$$

The function $\zeta(s,a)$ extends to a holomorphic function of $s$ on the whole complex plane with a simple pole at $s = 1$ (see [6]).

**Theorem 1**. Let

$$\zeta(s, a) = \sum_{n=0}^{\infty} \zeta_n(a) s^n, \quad |s| < 1, \qquad (2.2)$$

then

$$\zeta_n(a) = -1 + \sum_{k=n}^{\infty} (-1)^{k+1} \begin{bmatrix} k \\ n \end{bmatrix} \frac{B_{k+1}(a-1)}{(k+1)!}, \qquad (2.3)$$

where $B_n(a)$ are the Bernoulli polynomials and the series is semi-convergent in the sense of [8, p. 328].

When $a = 1$, $\zeta(s,1) = \zeta(s)$ is the Riemann zeta-function. Thus we have:

**Corollary**. If

$$\zeta(s) = \sum_{n=0}^{\infty} \zeta_n s^n, \quad |s| < 1, \qquad (2.4)$$

then

$$\zeta_n = -1 + \sum_{k=n}^{\infty} (-1)^{k+1} \begin{bmatrix} k \\ n \end{bmatrix} \frac{B_{k+1}}{(k+1)!}, \qquad (2.5)$$

where $B_n = B_n(0)$ are the Bernoulli numbers.

Note that

$$\zeta^{(n)}(0,a) = \zeta_n(a)n!, \qquad (2.6)$$

where the derivatives are for the variable $s$.

*Proof of the theorem.* We need two well-known facts ([6]):

$$\frac{e^{ax}}{e^x - 1} - \frac{1}{x} = \sum_{k=0}^{\infty} \frac{B_{k+1}(a)}{(k+1)!} x^k \quad (|x| < 2\pi) \qquad (2.7)$$

and

$$\zeta(-k, a) = \frac{-B_{k+1}(a)}{k+1}, \quad k = 0, 1, \ldots. \qquad (2.8)$$

Now let $a > 0$ be fixed. The residue of $\zeta(s,a)$ at $s = 1$ is $1$. Therefore, the function

$$f(x) = \zeta(-x, a) + \frac{1}{x+1} \qquad (2.9)$$

is entire. Set $f(x) = \sum_{n=0}^{\infty} a_n x^n$. Then according to (2.8)

$$f(n) = \zeta(-n, a) + \frac{1}{n+1} = \frac{-B_{n+1}(a)}{n+1} + \frac{1}{n+1}, \qquad (2.10)$$

and the ETF provides

$$\sum_{n=0}^{\infty} \frac{-B_{n+1}(a)}{(n+1)!} x^n + \sum_{n=0}^{\infty} \frac{x^n}{(n+1)!} = e^x \sum_{k=0}^{\infty} a_k \phi_k(x), \qquad (2.11)$$

which, in view of (2.7) can be written as

$$\frac{1}{x} - \frac{e^{ax}}{e^x - 1} + \frac{1}{x}(e^x - 1) = -e^x \left( \frac{e^{x(a-1)}}{e^x - 1} - \frac{1}{x} \right) = e^x \sum_{k=0}^{\infty} a_k \phi_k(x). \qquad (2.12)$$

The second equality, again in view of (2.10), turns into

$$\sum_{n=0}^{\infty} \frac{-B_{n+1}(a-1)}{(n+1)!} x^n = \sum_{k=0}^{\infty} a_k \phi_k(x). \qquad (2.13)$$

Substituting here (1.4) and comparing the coefficients in front of $x^k$ on both sides we arrive at the equations

$$\sum_{n=k}^{\infty} \{{}_k^n\} a_n = \frac{-B_{k+1}(a-1)}{(k+1)!}, \quad k = 0, 1, \ldots. \tag{2.14}$$

This is an infinite system for $a_n$ with a triangular matrix. For every $n = 0, 1, \ldots$, we multiply the $k$-th row ($\forall\, k \geq n$) by $(-1)^{k-n} [{}_n^k]$ and use the identity:

$$\sum (-1)^{k-n} [{}_n^k] \{{}_m^n\} = \delta_{k,m} \tag{2.15}$$

([see [7, p. 264]) to find

$$a_n = \sum_{k=n}^{\infty} (-1)^{k-n+1} [{}_n^k] \frac{B_{k+1}(a-1)}{(k+1)!}. \tag{2.16}$$

From the definition of $f(x)$ one has

$$\zeta(x, a) = \frac{-1}{1-x} + \sum_{n=0}^{\infty} (-1)^n a_n x^n \tag{2.17}$$

or, using the series expansion of $1/(1-x)$, $|x| < 1$,

$$\zeta(x, a) = \sum_{n=0}^{\infty} ((-1)^n a_n - 1) x^n. \tag{2.18}$$

Therefore,

$$\zeta_n(a) = -1 + (-1)^n a_n, \tag{2.19}$$

which combined with (2.16) leads to the desired result. The proof is completed.

In particular, when $n = 0$ one verifies that

$$\zeta(0, z) = \zeta_0(a) = -1 + (-1) B_1(a-1) = 1/2 - a, \tag{2.20}$$

as $B_1(x) = x - 1/2$.

When $n = 1$ we have $[{}_1^k] = (k-1)!$ and

$$\zeta_1(a) = -1 + \sum_{k=1}^{\infty} (-1)^{k+1} \frac{B_{k+1}(a-1)}{k(k+1)}. \tag{2.21}$$

At the same time (see [6])

$$\zeta_1(a) = \zeta'(0, a) = \log \Gamma(a) - \frac{1}{2} \log 2\pi \qquad (2.22)$$

which leads to the well-known representation [8, p.336]

$$\log \Gamma(1 + a) = \frac{1}{2} \log 2\pi - 1 + \sum_{k=1}^{\infty} (-1)^{k+1} \frac{B_{k+1}(a)}{k(k+1)}. \qquad (2.23)$$

Equation (2.23) comes, for instance, from the asymptotic representation

$$\log \Gamma(z + a) = (z + a - \frac{1}{2}) \log z - z + \frac{1}{2} \log(2\pi) + \sum_{k=1}^{\infty} (-1)^{k+1} \frac{B_{k+1}(a)}{k(k+1)} z^{-k}, \qquad (2.24)$$

(see [6, 1.18 (12)]) by setting $z = 1$.

When $n = 2$, we have $\left[ \begin{smallmatrix} k \\ 2 \end{smallmatrix} \right] = (k-1)! H_{k-1}$, where

$$H_{k-1} = 1 + \frac{1}{2} + \ldots + \frac{1}{k-1}, \qquad (2.25)$$

are the harmonic numbers. From the theorem

$$\zeta_2(a) = -1 + \sum_{k=2}^{\infty} (-1)^{k+1} H_{k-1} \frac{B_{k+1}(a-1)}{k(k+1)}, \qquad (2.26)$$

etc.

**Notes.** A representation of the coefficients $\zeta_n$ as certain limits is given in [3, p.215], [10]. For $\zeta''(0, a) = 2\zeta_2(a)$ see also the discussion on pp. 204-207 in [3]. Apostol [1] obtained a closed form of $\zeta_n$ in terms of Taylor's coefficients in the expansion of $\Gamma(s)\zeta(s) - 1/(s-1)$ about $s = 1$. Other computations of $\zeta_n$ can be found in [10]. The Taylor coefficients $\zeta_n$ are related to the Stieltjes constants $\gamma_n$ in the Laurent series of the Zeta function centered at $s = 1$ (see [9]).

## 3. The Lerch zeta function

The Lerch zeta function (or Lerch Transcendent) represents a generalization of the Hurwitz zeta function,

$$\Phi(\lambda, s, a) = \sum_{n=0}^{\infty} \frac{\lambda^n}{(n+a)^s} \ . \tag{3.1}$$

Here $|\lambda| \leq 1$ and $a > 0$. A detailed definition of $\Phi$ and its basic properties can be found in [6]. Assuming $\lambda \neq 1$, we show how Theorem 1 changes for this function. First we recall a class of functions $\beta_n(a, \lambda)$ introduced by Apostol [2] (see also [5]) and defined by the generating function

$$\frac{z e^{az}}{\lambda e^z - 1} = \sum_{n=0}^{\infty} \beta_n(a, \lambda) \frac{z^n}{n!} . \tag{3.2}$$

When $\lambda = 1$, $\beta_n(a, 1)$ are the Bernoulli polynomials. When $\lambda \neq 1$, $\beta_n(a, \lambda)$ are rational functions of $\lambda$ and polynomials in the variable $a$ of order $n-1$. Thus

$$\beta_0(a, \lambda) = 0, \ \beta_1(a, \lambda) = \frac{1}{\lambda - 1}, \ \beta_2(a, \lambda) = \frac{2a(\lambda - 1) - 2\lambda}{(\lambda - 1)^2}, \ldots \text{etc} . \tag{3.3}$$

The function $\Phi(\lambda, s, a)$ extends as a holomorphic function of $s$ on the entire complex plane. Apostol proved that for $s = -m$, $m = 0, 1, \ldots$,

$$\Phi(\lambda, -m, a) = \frac{-\beta_{m+1}(a, \lambda)}{m+1} , \tag{3.4}$$

which corresponds to (2.8).

Let $\lambda \neq 1$ and consider the Taylor series representation $\Phi(\lambda, s, a)$ in $s$

$$\Phi(\lambda, s, a) = \sum_{n=0}^{\infty} c_n(a, \lambda) s^n . \tag{3.5}$$

**Theorem 2**. The coefficients $c_n(a, \lambda)$ can be represented as semi-converegent series

$$c_n(a, \lambda) = \sum_{k=n}^{\infty} (-1)^{k-n+1} \begin{bmatrix} k \\ n \end{bmatrix} \frac{\beta_{k+1}(a-1, \lambda)}{(k+1)!} . \tag{3.6}$$

The proof follows the same steps as in Theorem 1. We apply the ETF (1.8) to the function $f(x) = \Phi(\lambda, -x, a)$ in order to obtain, in view of (3.4), the representation

$$-\sum_{n=0}^{\infty} \frac{\beta_{n+1}(a, \lambda)}{(n+1)!} x^n = e^x \sum_{k=0}^{\infty} c_k(a, \lambda) \phi_k(x). \tag{3.7}$$

Then since

$$e^{-x} \frac{e^{ax}}{\lambda e^x - 1} = \frac{e^{(a-1)x}}{\lambda e^x - 1}, \tag{3.8}$$

we find from (3.2) and (3.7)

$$-\sum_{n=0}^{\infty} \frac{\beta_{n+1}(a-1, \lambda)}{(n+1)!} x^n = \sum_{k=0}^{\infty} c_k(a, \lambda) \phi_k(x). \tag{3.9}$$

The rest of the proof follows by comparing coefficients for $x^k$ on both sides in (3.9)

$$\sum_{n=k}^{\infty} \begin{Bmatrix} n \\ k \end{Bmatrix} c_n(a, \lambda) = \frac{-\beta_{k+1}(a-1, \lambda)}{(k+1)!} , \quad k = 0, 1, \ldots \tag{3.10}$$

and solving this system for $c_n(a, \lambda)$ by using (2.15).